\documentstyle[12pt]{article}

\begin{document}
\newcommand{\ol }{\overline}
\newcommand{\ul }{\underline }
\newcommand{\ra }{\rightarrow }
\newcommand{\lra }{\longrightarrow }
\newcommand{\ga }{\gamma }
\newcommand{\st }{\stackrel }

\title{\Large\bf The Baer-invariant of a Semidirect Product}
\author{by\\ Behrooz Mashayekhy\\ Department of Mathematics \\
Ferdowsi University of Mashhad\\ P.O.Box 1159-91775\\ Mashhad, Iran\\e-mail:
mashaf@science2.um.ac.ir}
\date{ }

\maketitle

\begin{abstract}
 In 1972 K.I.Tahara [7,2 Theorem 2.2.5] , using cohomological method, showed
that if a finite group $G=T\rhd\!\!\! <N$ is the semidirect product of a normal
subgroup $N$ and a subgroup $T$ , then $M(T)$ is a direct factor of $M(G)$ ,
where $M(G)$ is the Schur-multiplicator of $G$ and in the finite case , is the
second cohomology group of $G$ . In 1977 W.Haebich [1 Theorem 1.7] gave another
proof using a different method for an arbitrary group $G$ .

 In this paper we generalize the above theorem . We will show that
${\cal N}_cM(T)$ is a direct factor of ${\cal N}_cM(G)$ , where
${\cal N}_c$ [3 page 102] is the variety of nilpotent groups of
class at most $c\geq 1$ and ${\cal N}_cM(G)$ is {\it the
Baer-invariant} of the group $G$ with respect to the variety
${\cal N}_c$ [3 page 107] .\\
A.M.S.Classification (1990): 20E22,20F12,20F18 
\end{abstract}

\begin{center}{\bf 1. Notations and Preliminaries} \end{center}
{\bf Definition 1.1}

  The group $G$ is said to be a {\it semidirect product} of a normal subgroup
$A$ and a subgroup $B$ , denoted by $G=B\rhd\!\!\! <A$ ( or a {\it splitting
extension} of $A$ by $B$ ) if

$(i)$ $G$ is generated by $A$ and $B$ ,

$(ii)$ $A\cap B=1$ .

 Since $A$ is normal in $G$ , the maps $\theta b:\  a\mapsto
a^b\ ,\ a\in A\ $ , fo all $b\in B$ are automorphisms of $A$ and they induces a
homomorphism $\theta :\ B\ra Aut(A)$ which is called the action of $B$ on $A$ .
$G$ is determined up to isomorphism by $\theta$ and is therefore called the
semidirect product of $A$ and $B$ under $\theta$ ( or the splitting extension
of $A$ by $B$ under $\theta$ ) . Note that every element of $G$ is uniquely
determined by $ab$ , for $a\in A$ and $b\in B$ .

 Now the following elementary results on semidirect product are needed in our
work , see [1 page 421, 5] \\
{\bf Lemma 1.2}

 Let $G$ be the semidirect product of $A$ and $B$ under $\theta :B\ra Aut(A)$
and $\ol G$ be the semidirect product of $\ol A$ and $\ol B$ under $\ol \theta$
. If
$$ \alpha :A\lra {\ol A} \ \ \ \  {\rm and} \ \ \ \ \beta :B\lra {\ol B} $$
are epimorphisms such that
$$ \alpha((\theta b)(a))=({\ol \theta}\beta b)(\alpha a)\ ,\ {\rm for\ all}\
a\in A\ ,\ {\rm and}\ b\in B\ ,$$
then the map
$$ \tau :G\lra {\ol G} $$
$$ \ \ \ \ \ \ \ \ \  ab\longmapsto (\alpha a)(\beta b) $$
is an epimorphism extending $\alpha$ and $\beta$ .\\
\ \ \\  \ \ \\
{\bf Lemma 1.3}

 Let $G$ be the semidirect product of $A$ and $B$ under $\theta :B\ra Aut(A)$ .
If $N$ is a subgroup of $A$ which is normal in $G$ , then $G/N$ is the
semidirect product of $A/N$ and $B$ under
$$ \tilde{\theta}:B\lra Aut(A/N)$$
$$ b\longmapsto \tilde{\theta}b \ \ \ \ \  $$
where $\tilde{\theta}b:A/N\lra A/N$ given by $aN\mapsto \theta(a)N$ .\\

 In the following theorem a free presentation for $G$ is introduced in terms of
free presentations of $A$ and $B$ .\\
{\bf Theorem 1.4}

 Let $G$ be a semidirect product of $A$ and $B$ under $\theta :B\lra Aut(A)$ and
$$ 1\lra R_1\lra F_1\st{\nu_1}{\lra}A\lra 1\ \ \ ,\ \ \ 1\lra R_2\lra
F_2\st{\nu_2}{\lra}B\lra 1 $$
be free presentations for $A$ and $B$ , respectively. Then
$$ 1\lra R\lra F\lra G\lra 1 $$
is a free presentation for $G$ , where

$(i)$ $F=F_1*F_2$ , the free product of $F_1$ and $F_2$ ;

$(ii)\ R=R_1^FR_2^FS\ ;$

$(iii)\ S=<f_1^{-1}\ol{f_1}[f_2,f_1]\ |\ f_1,\ol{f_1}\in F_1;f_2\in
F_2;\nu_1\ol{f_1}=\theta (\nu_2f_2)(\nu_1f_1)>^F\ .$\\
{\bf Proof.}

 See[1 Lemma 1.4] . $\Box$  \\

 Let $\cal V$ be a variety of groups defined by a set of laws $V$ and $G$ be a
group with a free presentation
$$ 1\lra R\lra F\lra G\lra 1\ \ .$$
Then  {\it the Baer-invariant} of $G$ with respect to the variety $\cal V$,
denoted by ${\cal V}M(G)$ , is defined to be
$$ \frac{R\cap V(F)}{[RV^*F]} $$
where $V(F)$ is the verbal subgroup of $F$ and $[RV^*F]$ is the least normal
subgroup of $F$ , contained in $R$ , generated by
$$\{v(f_1,\ldots ,f_ir,\ldots ,f_s)v(f_1,\ldots ,f_s)^{-1}|r\in R,f_i\in F,v\in
V,1\leq i\leq s\}\ .$$
It is easily seen that the Baer-invariant of the group $G$ with respect to the
variety $\cal V$ is always abelian and that it is independent of the choice of
the free presentation of $G$ , see [3 Lemma 1.8].\\
In particular, if $\cal V$ is the variety of abelian groups, then the
Baer-invariant of $G$ with respect to $\cal V$ will be
$$ \frac{R\cap F'}{[R,F]}\ \ ,$$

which , by I.Schur , is isomorphic to the Schur-multiplicator of $G$ , denoted by
$M(G)$ , in general , and in the finite case , is isomorphic to the second
cohomology group of $G$ , $H^2(G,{\bf C})$ [2 Theorem 2.4.6].

 If $\cal V$ is the variety ${\cal N}_c$ of nilpotent groups of
class at most $c\geq 1$ , then the Baer-invariant of $G$ with respect to ${\cal
N}_c$ is
$$ \frac{R\cap \ga_{c+1}(F)}{[R, _cF]}\ \ ,$$
where $\ga_{c+1}(F)$ is the $(c+1)$st-term of the lower centeral series and
$[R, _cF]$ stands for $[R,\underbrace {F,F,\ldots ,F}_{c-times}]$ .
For further details , properties , conventions , see [3,4,5] .\\
\begin{center} {\bf 2. The Main Result}  \end{center}

 The following theorem is fundemental to the proof of Theorem 2.2 . We adopt the
notations and conventions from section 1 , in what follows.\\
\ \ \\ \ \ \\
{\bf Theorem 2.1}

$(i)\ R_1$ and $[R_2,F_1]$ are subgroups of $S$ ;

$(ii)\ R=R_2S$ ;

$(iii)\ R\cap \ga_{c+1}(F)=(R_2\cap \ga_{c+1}(F_2))(S\cap \ga_{c+1}(F))$ , for
all $c\geq 1$ ;

$(iv)\ [R,\ _cF]=[R_2,\ _cF_2]\prod  [R_2,F_1,F_2]_c[S,\ _cF]$ , for
all $c\geq 1$ , \\
where
$$ \prod [R_2,F_1,F_2]_c=<[r_2,f_1,\ldots ,f_c]\ |\ f_i\in F_1\cup
F_2,r_2\in R_2,1\leq i\leq c,\exists k,\ f_k\in F_1>^F\ \ .$$
In particular, $\prod  [R_2,F_1,F_2]_1=[R_2,F_1]$ .\\
{\bf Proof.}

$(i)$ If $r_1\in R_1$ , then $\nu_1r_1=1$ and hence
$\theta (\nu_2f_2)(\nu_1f_1)=1$ , for all $f_2\in F_2$.  Therefore
$r_1=1^{-1}r_1[f_2,1]$ , so $r_1\in S$ i.e $R_1\leq S$ .\\
If $r_2\in R_2$ ,then $\nu_2r_2=1$ , so $\theta (\nu_2r_2)$ is the identity
automorphism. Thus $[r_2,f_1]=f_1^{-1}f_1[r_2,f_1]\in S$ and so $[R_2,F_1]\leq
S$ .

$(ii)$
$$ \ \ \ \ \ \ \ R=R_2^FR_1^FS\ \ \ \ \ \ , \ \ \ {\it by\ Theorem\ 1.4} $$
$$\ \ =R_2^FS\ \ \ \ \ \ \ \ \ \ \ , \ \ \ \ \ \ \ {\it by\ (i)} $$
$$\ \ \ \ =R_2[R_2,F]^FS  \ \ \ \ \ \ \ \ \ \ \ \ \ \ \ \ \ \ \ \ \ $$
$$\ \ \ \ \ \ \ \ \ \ =R_2[R_2,F_1]^FS\ \ \ \ , \ \ \ {\it since\ } R_2\unlhd
F_2 $$ $$ =R_2S\ \ \ \ \ \ \ \ , \ \ \ {\it by\ (i)}\ \ .$$

$(iii)$ Since $F=F_1*F_2$ , we have
$$ \ga_{c+1}(F)=\ga_{c+1}(F_1)\ga_{c+1}(F_2)\prod [F_1,F_2]_{c+1}\ \ ,$$
where
$$ \prod [F_1,F_2]_{c+1}=<[F_1,F_2,F_{i_1},\ldots ,F_{i_{c-1}}]\ |\ i_j\in
\{1,2\},1\leq j\leq c-1>$$
and $\prod [F_1,F_2]_{c+1}\unlhd F$ ( to find a proof see
M.R.R.Moghaddam
[4] ). Also we know that $F=F_2\rhd\!\!\!< F_1[F_1,F_2]$ is the semidirect
product of $F_2$ and $F_1[F_1,F_2]$ ( since $F_2\cap F_1[F_1,F_2]=1$ and
$F_1[F_1,F_2]\unlhd F$ ) and $S\leq F_1[F_1,F_2]$. So using part $(ii)$ and
the above remarks, we have
$$ R\cap \ga_{c+1}(F)=R_2S\cap \ga_{c+1}(F_1)\ga_{c+1}(F_2)\prod
[F_1,F_2]_{c+1}\ \ \ \ \  $$
$$\ \ \ \ \ \ \ \ \ \ \ \ \ \ \ \ =(R_2\cap \ga_{c+1}(F_2))(S\cap \ga_{c+1}(F_1)
\prod [F_1,F_2]_{c+1}) $$
$$ \ \ =(R_2\cap \ga_{c+1}(F_2))(S\cap \ga_{c+1}(F))\ \ \ \ \ $$

$(iv)$ Use induction on $c$ . Let $c=1$ . Then
$$ [R,F]=[R_2S,F]\ \ \ \ , \ \ \ \ {\it by\ (ii)} $$
$$\ \ \ \ \ \ \ \ \ \ \ \ \ \ =[R_2,F][S,F]\ \ \ , \ \ {\it since\ } S\unlhd F
$$
$$\ \ \ \ \ \ \ \ \ \ \ \ \ \ \ \subseteq [R_2,F_2][R_2,F_1]^F[S,F]\ \ , \ {\it
since\ } F=F_1*F_2 $$
$$ \ \ \ \ \ \ \ \ \ \ \ \ \ \ \ \ \ \ =[R_2,F_2][R_2,F_1][S,F]\ \ \ , \ {\it
since\ } [R_2,F_1]\leq
S\ \ . $$ Clearly \ \ \ \ $[R_2,F_2][R_2,F_1][S,F]\subseteq [R,F]$\ \ .
Hence\ \ $[R,F]=[R_2,F_2][R_2,F_1][S,F]$\ \ .

 Now, suppose $[R,\ _kF]=[R_2,\ _kF_2]\prod [R_2,F_1,F_2]_k[S,\ _kF]$ .
 Then we have
$$ [R,\ _{k+1}F]=[[R,\ _kF],F]$$    $$=[[R_2,\ _kF_2]\prod [R_2,F_1,F_2]_k[S,\
_kF],F]\ \ \ \ $$
$$\ =[[R_2,\ _kF_2],F][\prod [R_2,F_1,F_2]_k,F][[S,\ _kF],F]\ \ \ ,\ \ \
(by\ induction \ hypothesis) $$
$$\ \ \ \ \ \subseteq [[R_2,\ _kF_2],F_2]\prod [R_2,F_1,F_2]_{k+1}[S,\ _{k+1}F]
 \ \ \ \ , \ \ \ \ (since\ [S,\ _kF]\ ,\ \prod [R_2,F_1,F_2]_k \unlhd F) $$
$$\ \ \ \ \ \ \ \ \ \ \subseteq [R,\ _{k+1}F]\ \ \ .$$
Therefore , by induction we have
$$ [R,\ _cF]=[R_2,\ _cF_2]\prod [R_2,F_1,F_2]_c[S,\ _cF]\ \ \ \ for\ all\ c\geq
1\ .\ \ \Box$$\\

Now we are in a position to state and prove the main theorem of this paper.\\
{\bf Theorem 2.2}

 Let $G$ be a semidirect product of $A$ by $B$
under $\theta :B\lra Aut(A)$ (or a splitting extension of $A$ by $B$ under
$\theta $) , and ${\cal N}_c$ be the variety of nilpotent
groups of class at most $c\ (c\geq 1)$. Then
$$ {\cal N}_cM(G)\cong {\cal N}_cM(B)\oplus \frac {S\cap \ga_{c+1}(F)}{\prod
[R_2,F_1,F_2]_c[S,\ _cF]}\ \ .$$
In particular , ${\cal N}_cM(B)$ can be regarded as a direct factor of
${\cal N}_cM(G)$ .\\
{\bf Proof.}

 By the previous assumptions and notations we have
$$F\st{\varphi}{\lra} \frac{F}{[R_2,\ _cF_2]^F}\st{\eta}{\lra}\frac{F}{[R_2,\
_cF_2]\prod [R_2,F_1,F_2]_c[S,\ _cF]}\ \ , $$
where $\varphi$ and $\eta$ are natural homomorphisms. Then for any $c\geq 1$ ,
we have $$ \frac{R\cap \ga_{c+1}(F)}{[R,\ _cF]}\cong (\eta \varphi )(R\cap
\ga_{c+1}(F)) $$
$$ \cong (\eta \varphi )(R_2\cap \ga_{c+1}(F_2))(\eta
\varphi )(S\cap \ga_{c+1}(F))\ \ , \  by\ Theorem\ 2.1\ (iii)\ \ \ \ \ \
(*)\ \ \ . $$ Cosider the following two natural homomorphisms
$$
\frac{F_1*F_2}{[R_2,\ _cF_2]^F} \st{h}{\lra}  F_1*\frac{F_2}{[R_2,\ _cF_2]}
\st{g}{\lra} \frac{F_1*F_2}{[R_2,\ _cF_2]^F}\ \ \ \ \ , $$
given by
$$\ol{f_1}  \longmapsto   f_1\ \ \ \ \ \ \ \ \ \ f_1  \longmapsto  \ol{f_1}
\ \ \ \ $$
$$ \ol{f_2} \longmapsto  \ol{f_2}\ \ \ \ \ \ \ \ \ \ \ol{f_2} \longmapsto
\ol{f_2} \ \ \ . $$
Clearly  \ \ \ \ $h\circ g=1\ \ \ \&\ \ \ g\circ h=1 $ i.e $g$ is the inverse of
$h$ and so $h$ is an isomorphism. Thus we have
$$ \frac{F_1*F_2}{[R_2,\ _cF_2]^F}=\varphi (F)\cong F_1*\frac{F_2}{[R_2,\
_cF_2]}\ \ \ .$$
Also $\varphi (F_2)={F_2}/{[R_2,\ _cF_2]}$ and $$\varphi (F_1[F_1,F_2])\cong
\varphi (F_1)[\varphi (F_1),\varphi (F_2)]\cong F_1[F_1,F_2/[R_2,\ _cF_2]]\ \ .
$$
Therefore
$$ \varphi (F)\cong F_1*\frac{F_2}{[R_2,\ _cF_2]}\cong \frac{F_2}{[R_2,\
_cF_2]}\rhd\!\!\! <F_1[F_1,\frac{F_2}{[R_2,\ _cF_2]}] $$
 $$ \ \ \ \ \ \ \ \ \ \ \ \ \ \ \ \ \ \ \ \ \ \ \ \ \ \
\cong \varphi(F_2)\rhd\!\!\!<\varphi(F_1)[\varphi(F_1),\varphi(F_2)]\ \ .\ \ \ \ \ \
\ \ \ \ \ \ \ \ \ \ \ \ \ \ \ \ \  (**)$$
Thus by Lemma 1.3 and the property that $Ker(\eta )\leq
\varphi(F_1)[\varphi(F_1),\varphi(F_2)]\ $ we have
$$ (\eta \varphi )(F)\cong \frac{\varphi(F)}{Ker(\eta )}=\frac{\varphi(F)}{
\varphi(\prod [R_2,F_1,F_2]_c[S,\ _cF])}\cong
\varphi(F_2)\rhd\!\!\! <\frac{\varphi(F_1)[\varphi(F_1),\varphi(F_2)]}{Ker(\eta
)}\ \ , $$
Clearly $(\eta \varphi )(F_2)\cong \varphi(F_2)$ and $(\eta\varphi )(F_1)\cong
\varphi(F_1)/Ker(\eta  )$ , thus we have
$$ (\eta\varphi )(F)\cong
\varphi(F_2)\rhd\!\!\! <\frac{\varphi(F_1)[\varphi(F_1),\varphi(F_2)]}{Ker(\eta
)}\ \ \ \ (by\ 1.3\ and\ (**)) $$
$$\ \ \ \ \ \ \ \ \ \ \ \cong (\eta\varphi )(F_2)\rhd
\!\!\! <(\eta\varphi )(F_1)[(\eta\varphi )(F_1),(\eta\varphi )(F_2)]\ \ \ .$$
So
$$(\eta\varphi )(R_2\cap \ga_{c+1}(F_2))\ \cap\ (\eta\varphi )(S\cap
\ga_{c+1}(F))\subseteq (\eta\varphi )(F_2)\cap
(\eta\varphi )(F_1)[(\eta\varphi )(F_1),(\eta\varphi )(F_2)]=1 $$
Hence, by $(*)$
$$ \frac{R\cap \ga_{c+1}(F)}{[R,\ _cF]}\cong (\eta\varphi )(R_2\cap
\ga_{c+1}(F_2))\oplus(\eta\varphi )(S\cap \ga_{c+1}(F)) $$
and
$$ (\eta\varphi )(R_2\cap\ga_{c+1}(F_2))=\frac{(R_2\cap
\ga_{c+1}(F_2))Ker(\eta\varphi )}{Ker(\eta\varphi )}\cong \frac{R_2\cap
\ga_{c+1}(F_2)}{(R_2\cap \ga_{c+1}(F_2))\cap Ker(\eta\varphi )} $$
$$\ \ \ \ \ \ \ \ \cong \frac{R_2\cap \ga_{c+1}(F_2)}{[R_2,\ _cF_2]}\cong {\cal
N}_cM(B)\ \ .$$
Also
$$ (\eta\varphi )(S\cap \ga_{c+1}(F))=\frac{(S\cap
\ga_{c+1}(F))Ker(\eta\varphi )}{Ker(\eta\varphi )}\cong \frac{S\cap
\ga_{c+1}(F)}{(S\cap \ga_{c+1}(F))\cap Ker(\eta\varphi )} $$
$$\ \ \ \ \ \ =\frac{S\cap \ga_{c+1}(F)}{\prod [R_2,F_1,F_2][S,\ _cF]}\ \ \ .$$
Therefore
$$ {\cal N}_cM(B\rhd\!\!\! <A)\cong \frac{R\cap \ga_{c+1}(F()}{[R,\ _cF]}\cong
{\cal N}_cM(B)\oplus \frac{S\cap \ga_{c+1}(F)}{\prod [R_2,F_1,F_2]_c[S,\ _cF]}\
\ .\ \ \ \Box$$

 Now we obtain the following corollaries:\\
{\bf Corollary 2.3} (Tahara [7,2 Theorem 2.2.5])

 Let $G=B\rhd\!\!\! <A$ be the semidirect product of a normal subgroup $A$ and
a subgroup $B$ . Then $M(B)$ is a direct factor of $M(G)$ .\\
{\bf Corollary 2.4} (Haebich [1 Theorem 1.7])

 Suppose $G=B\rhd\!\!\! <A$ is a semidirect product of $A$ by $B$ under
$\theta :B\ra Aut(A) $. By the notation of Theorem 1.4 we have
$$ M(G)\cong M(B)\oplus \frac{S^F\cap F'}{[R_2,F_1][S,F]} \ \ \ .  $$


\newpage
\begin{thebibliography}{7}

\bibitem{1} W.Haebich;``The Multiplicator of a Splitting Extension.''J.Algebra,
{\bf 44},(1977),420-33.

\bibitem{2} G.Karpilovsky;``The Schur Multiplier.'' London Math. Soc.Monographs
(New Series No.{\bf 2})(1987).

\bibitem{3} C.R.Leedham\_Green and S.McKay;``Baer-invariant,Isologism,Varietal
Laws and Homology.'' Acta Math.,{\bf 137}(1976),99-150.

\bibitem{4} M.R.R.Moghaddam;``The Baer-invariant of a Direct Product.''Archiv.
der. Math.,vol.{\bf 33}(1979),504-511.

\bibitem{5} D.J.S.Robinson;``A Course in the Theory of Groups.'' Springer Verlag,
(1982).

\bibitem{6} I.Schur;``Untersuchungen \"{u}ber die Darstellung der Endlichen Gruppen
durch Gebrochenen Linearen Substitutionen.''J.Reine Angew. Math.,{\bf 132}(1907)
,85\_137.

\bibitem{7} K.I.Tahara;``On the Second Cohomology of Semidirect Product.'' Math.
Z.,{\bf 129}(1972),365\_379.

\end{thebibliography}
\end{document}